\documentclass[leqno]{siamltex}
\usepackage{amsmath}
\usepackage{graphicx}
\usepackage{mathrsfs}
\usepackage{float}
\usepackage{amsfonts,amssymb,,graphicx}
\usepackage{dsfont}
\usepackage{pifont}
\usepackage{multirow}
\usepackage{float}
\usepackage{diagbox}
\usepackage{multicol}
\usepackage{arydshln}
\usepackage{bbm}
\usepackage[ruled]{algorithm2e}
\usepackage{booktabs}
\usepackage{hyperref}
\hypersetup{hypertex=true,
	colorlinks=true,
	linkcolor=magenta,
	anchorcolor=blue,
	citecolor=blue}
\usepackage{cite}
\usepackage{tikz}
\usepackage{hyperref}
\usepackage{mathtools}
\usepackage{mfirstuc}

\usepackage{subfigure}
\usepackage{color}

\def\3bar{{|\!|\!|}}

\usepackage{caption}

\title{Deep Finite Volume Method for Partial Differential Equations
}
\author{Jianhuan Cen\footnotemark[1]
		\and Qingsong Zou\footnotemark[2]
		}
\date{}
\begin{document}
\maketitle
\renewcommand{\thefootnote}{\fnsymbol{footnote}}
\footnotetext[1]{School of Computer Science and Engineering, Sun Yat-sen University, Guangzhou, 510006, China.}

\footnotetext[2]{Corresponding author. School of Computer Science and Engineering, and Guangdong Province Key Laboratory of Computational Science, Sun Yat-sen University, Guangzhou 510006, China. Email: mcszqs@mail.sysu.edu.cn.}

\captionsetup[figure]{labelfont={bf},labelformat={default},labelsep=period,name={Figure}}
\captionsetup[table]{labelfont={bf},labelformat={default},labelsep=period,name={Table}}

\begin{abstract}
In this paper, we introduce the Deep Finite Volume Method (DFVM), an innovative deep learning framework tailored for solving high-order (order \(\geq 2\)) partial differential equations (PDEs). Our approach centers on a novel loss function crafted from local conservation laws derived from the original PDE, distinguishing DFVM from traditional deep learning methods. By formulating DFVM in the weak form of the PDE rather than the strong form, we enhance accuracy, particularly beneficial for PDEs with less smooth solutions compared to strong-form-based methods like Physics-Informed Neural Networks (PINNs). A key technique of DFVM lies in its transformation of all second-order or higher derivatives of neural networks into first-order derivatives which can be comupted directly using Automatic Differentiation (AD). This adaptation significantly reduces computational overhead, particularly advantageous for solving high-dimensional PDEs. Numerical experiments demonstrate that DFVM achieves equal or superior solution accuracy compared to existing deep learning methods such as PINN, Deep Ritz Method (DRM), and Weak Adversarial Networks (WAN), while drastically reducing computational costs. Notably, for PDEs with nonsmooth solutions, DFVM yields approximate solutions with relative errors up to two orders of magnitude lower than those obtained by PINN. The implementation of DFVM is available on GitHub at \href{https://github.com/Sysuzqs/DFVM}{https://github.com/Sysuzqs/DFVM}.
\end{abstract}



\pagestyle{myheadings}\

\section{Introduction}\
Partial differential equations (PDEs) are prevalent and extensively applied in science, engineering, economics, and finance. Traditional numerical methods, such as the finite difference method \cite{Le2007}, the finite element method \cite{ZTZ2005}, and the finite volume method \cite{XZ2009}, have achieved great success in solving PDEs. However, traditional methods face significant challenges when dealing with certain nonlinear or high-dimensional PDEs, or complex domain problems. For instance, traditional numerical methods often suffer from the so-called ``curse of dimensionality," wherein the number of unknowns grows exponentially with the increase in dimension. Recently, deep learning methods have gained considerable popularity in solving PDEs due to their simplicity and flexibility \cite{DGM,DRM,ZZS2019, LJP2021,LKA2020}.

The loss functions in deep learning methods for PDEs can be broadly categorized into two types. The first type is directly designed according to the strong form of the PDE, and it is based on the residuals of the PDE in the least squares sense. A well-known example of this is the loss function used in Physics-Informed Neural Networks (PINNs) \cite{PINNs}. When the solution of the PDE is sufficiently smooth, deep solvers based on this type of loss function usually achieve high accuracy. Due to their simplicity, straightforward nature, and elegance, deep solvers employing this type of loss function have been widely applied to various problems, including the Schrödinger equation \cite{HZE2021}, the Hamilton-Jacobi-Bellman equations \cite{EHJ2017}, the Black-Scholes equation \cite{GC2019}, and problems with random uncertainties \cite{ZLG2019,ZZS2019}.

The second type of loss is designed according to the weak forms of the original PDE. This approach recognizes that the actual solution of a real physical problem might not be sufficiently smooth to satisfy the strong form of the PDE strictly, but it can satisfy a certain weak form of the PDE. Many weak-form-based losses have been proposed in the literature. For example, the Deep Ritz Method (DRM) \cite{DRM} employs an energy functional as its loss function. Weak Adversarial Networks (WAN) \cite{WAN} and variational PINNs (vPINN) \cite{vPINNs} both utilize loss based on the Petrov-Galerkin framework. WAN utilizes two neural networks to fit trial and test functions separately, whereas vPINN employs a single neural network to fit trial functions and utilizes polynomials for test functions. The loss of the weak PINNs (wPINN) \cite{wPINNs} is based on the well-known family of Kruzkhov entropies so that the wPINN approximates the entropy solution of the original problem.  The Deep Mixed Residual Method (MIM) \cite{MIM2022} uses the loss based on a lower-order PDE system derived from the original PDE.

Since a solution that satisfies the strong form of the PDE must also satisfy its weak form, but a solution that satisfies the weak form of the PDE may not satisfy the strong form, the weak-form type loss theoretically has a wider range of applications than the strong-form type loss. Many numerical experiments demonstrate that when the solution of a PDE is not sufficiently smooth, as is common in many physical scenarios, weak-form-based deep solvers have distinct advantages. For example, the wPINN can handle shocks and discontinuities effectively, ensuring a more accurate and stable solution even when the solution lacks smoothness.


In this paper, we propose a novel weak form loss. To illustrate our basic idea of loss design, we take the Poisson equation $ -\Delta u  = f$ as an example. When the practical solution $u$ is not in $C^2$ in the whole domain $\Omega$, the strong form equation 
\begin{equation}\label{strong}
 (\Delta u +f)({\bf x}) = 0 
\end{equation}
does not hold for all points ${\bf x}$ in $\Omega$. However, due to the local conservation of the flux, the solution $u$ satisfies the weak form  
\begin{equation}\label{weak1}
\int_{\partial V} \frac{\partial u}{\partial {\bf n}} {\text d}s +\int_V f {\text d} {\bf x}=0
\end{equation}
for all subdomain $V\subset \Omega$, where $\partial V$ is the boundary of $V$ and ${\bf n}$ is the unit outward normal direction. Based on this observation, we design the novel  loss as 
\begin{equation}\label{lossnovel}
{\mathcal J} (u) =\sum_{{\bf x}_0\in {\mathcal S}} \left|\int_{\partial V_{{\bf x}_0}} \frac{\partial u}{\partial {\bf n}} {\text d}s +\int_{V_{{\bf x}_0}} f {\text d} {\bf x}\right|^2.
\end{equation}
Here ${\mathcal S}$ is a certian set of points sampling from $\Omega$, $V_{\bf x_0}\subset \Omega$ is a so-called {\it control volume} surrounding ${\bf x_0}$.

Since the loss \eqref{lossnovel} depends on a finite number of volumes 
$V_{\bf x_0}, {\bf x}_0\in {\cal S}$, we refer to the deep learning method based on the loss \eqref{lossnovel} as the Deep Finite Volume Method (DFVM). As will be explained later, the DFVM can be extended from solving the above Poisson equation to general second-order PDEs and even to any higher-order PDEs, as higher-order PDEs can be transformed into second-order systems by introducing intermediate variables.

Let's now introduce the significance of the DFVM. 
Firstly, as a deep solver based on a weak form loss, the DFVM achieves higher accuracy  than methods using strong-form loss in solving many real-world physical problems, especially those involving singularities or discontinuities.  Its alignment with the law of conservation ensures consistent and accurate performance across a broad range of applications, making it particularly effective in complex scenarios.

Secondly,
Compared to other weak form-based deep solvers, the DFVM offers distinct advantages through three key innovations:
1)Unlike methods that simulate both trial and test functions separately, such as WAN, V-PINN, and wPINN, which employ different strategies involving neural networks and polynomials, the DFVM simplifies by using characteristic functions associated with control volumes. These volumes, geometrically shaped like cubes or balls centered at points, replace the need for complex test function simulations, thereby enhancing computational efficiency.2) The DFVM achieves higher solution accuracy by effectively capturing local conservation properties. Extensive numerical experiments demonstrate its superiority in accuracy over other weak-based methods across both singular and non-singular cases. 3) While methods like DRM and wPINNs face limitations with asymmetric PDEs or specific applications, the DFVM excels in solving a broad range of PDEs. It handles asymmetric, higher-order, and high-dimensional equations effectively, making it a versatile tool for solving complex PDEs across various domains. These innovations position the DFVM as a robust and efficient approach in the realm of weak form-based deep solvers, emphasizing computational simplicity, high accuracy, and broad applicability across diverse physical problems.

The third important advantage of the DFVM is its efficiency in solving high-dimensional PDEs compared to the popular PINN, even when the true solution is very smooth. To illustrate why this is the case, let's provide some details on the practical calculation of the PINN loss and the DFVM loss. Taking again the Poisson equation \eqref{strong} as an example, the loss of the PINN involves computing  $\Delta u_{\theta}$, where $u_\theta$ is the neural network to approximate the solution $u$. Recalling that in deep learning-related calculus, the first-order derivative of a neural network function is often implemented using the Automatic Differentiation mechanism (AD), which is very effective and powerful because the calculation of first-order derivatives only requires one reverse-pass operation after the forward-pass operation for calculating the function value has been completed. However, calculating second-order derivatives is almost equivalent to applying first-order AD $d+1$ times for a d-dimensional neural network function. This results in the computational cost increasing almost proportionally with the dimension $d$.
To justify our above observation,  we tested a function $u_\theta:\mathbb{R}^d \to \mathbb{R}$, represented by a fully connected network with 6 hidden layers and 200 neurons per layer. Table \ref{table:AD} lists the computing times for evaluating $u_\theta, \nabla u_\theta$, and $\Delta u_\theta$ 
1000 times for each of 1000 random sampling points in $\mathbb{R}^d$  using a computer equipped with an NVIDIA TITAN RTX. We observe that the computational cost of first-order derivatives is relatively consistent across different dimensions, whereas the cost of $\Delta u_\theta$ increases significantly with the dimension.
Since training a deep PDE solver often involves many calculations of the loss and its derivative (with respect to weights) on a large number of sampling points, this indicates that for solving high-dimensional PDEs, one should try to avoid calculating second or higher-order derivatives of a network directly using AD.
\begin{table}[H]
    \caption{Computing time by the AD.}
    \label{table:AD}
    \centering
    \begin{tabular}{rccc} 
     dim & $u_\theta$ & $\nabla u_\theta$ &  $\Delta u_\theta$ \\ [0.5ex] 
     \hline
    2   & 0.52s & 0.41s  & 2.32s \\
    10  & 0.58s & 0.42s  & 10.81s \\ 
    50  & 0.67s & 0.44s  & 53.66s \\
    100 & 0.75s & 0.50s  & 110.20s \\
    200 & 0.77s & 0.54s  & 219.40s \\ 
     \hline
    \end{tabular}
\end{table}

The calculation of the DFVM  loss \eqref{lossnovel} involves 
computing the integral $\int_{\partial V_{{\bf x}_0}} \frac{\partial u}{\partial {\bf n}} {\text d}s$. 
In practice, the integral will be computed with a certain numerical quadrature such as   
\begin{equation}\label{disLaplace}
   \int_{\partial V} \frac{\partial u_\theta}{\partial {\bf n}} ds \approx Q(\frac{\partial u_\theta}{\partial {\bf n}},\partial V) :=\sum_{{\bf y}_j\in {\mathcal N}} c_j \frac{\partial u_\theta}{\partial {\bf n}}({\bf y}_j), 
 \end{equation}
where $c_j, {\bf y}_j$ are weights and locations of the selected quadrature,  the control volume $V_{\bf x_0}$ is often selected as a cube or ball centered as ${\bf x}_0$, ${\mathcal N}\subset \partial V$ is the set of the nodes in the quadrature. Next we analyze the computational cost for \eqref{disLaplace}. We observe first that the formula in \eqref{disLaplace} 
only involves computing first order derivatives of $u_\theta$.  Table \ref{table:AD} shows that the cost for a single first order derivative term is independent of the dimension $d$.  Then the total cost for the quadrature in \eqref{disLaplace} depends mainly upon $\#{\mathcal N}$, the number of nodes in the quadrature. Considering that the size of the volume $V_{\bf x}$ is often relatively small (i.e. ${\mathcal O(10^{-5})}$), we may choose a quadrature such that $\#{\mathcal N}$ is 
not too large to achieve a very accurate approximation of $\int_{\partial V} \frac{\partial u_\theta}{\partial {\bf n}} ds$. Namely in the case that the dimension $d$ is very large, the cost for computing for  $\int_{\partial V} \frac{\partial u_\theta}{\partial {\bf n}} ds$ might be far less than that $\Delta u_\theta$.  Consequently, the DFVM might be more efficient than the PINN in handling high-dimensional PDEs.

We may explore the practical DFVM loss in the following different perspective.  The formula \eqref{disLaplace} actually introduces a novel method for calculating the Laplace operator of a neural network. That is, for a sufficiently small  $h$ and for an arbitrary point ${\bf x}_0\in \Omega$,
\begin{equation}\label{disLaplace1}
\Delta u_\theta({\bf x}_0)\approx \frac{1}{|V|} \int_{V} \Delta u_{\theta}({\bf x}) d{\bf x} = \frac{1}{|V|} \int_{\partial V} \frac{\partial u_\theta}{\partial {\bf n}} ds \approx   \frac{1}{|V|}\sum_{{\bf y}_j\in {\mathcal N}} c_j \frac{\partial u_\theta}{\partial {\bf n}}({\bf y}_j), 
\end{equation}
where $V=V_{{\bf x}_0, h}$. 
Scheme \eqref{disLaplace1} approximates \(\Delta u_\theta(\mathbf{x}_0)\) using first-order derivative terms computed via AD, differing from both direct AD calculations and traditional finite difference schemes. 
With the scheme \eqref{disLaplace1}, we can update the PINN loss to obtain a novel loss
\begin{equation}\label{lossnovel1}
{\mathcal J} (u_\theta) = \sum_{{\bf x}_0 \in {\mathcal S}} \left| \frac{1}{|V|} \sum_{{\bf y}_j \in {\mathcal N}} c_j \frac{\partial u_\theta}{\partial {\bf n}}({\bf y}_j) + f ({\bf x}_0) \right|^2,
\end{equation}
which is actually a special case of our DFVM loss \eqref{lossnovel} in which we use  
the numerical scheme
\[
\int_{V_{\bf x_0}} f({\bf x}) d{\bf x} \approx f({\bf x_0}) |V_{\bf x_0}|.
\]
 It is worth noting that schemes  similar to \eqref{disLaplace1} and \eqref{lossnovel1}  can be designed also for general second-order differential equations or even higher-order equations.

We conduct numerous numerical experiments to validate the effectiveness of DFVM. Firstly, in Section 3.1, we assess the efficiency of \eqref{disLaplace1} on randomly chosen networks \( u_\theta \). We will  show that with an appropriate \( h \) (e.g., \( 10^{-5} \)), it achieves a mean absolute error (MAE) of approximately \( 10^{-11} \) compared to AD, while requiring less computational time.  Secondly, we apply DFVM to solve various types of PDEs, such as the Poisson equation across four distinct cases, highlighting its capabilities in handling singular problems, high-dimensional equations, and complex domains, including adaptive strategies. Additionally, we demonstrate DFVM's efficacy in solving higher-order equations like the biharmonic and Cahn-Hilliard equations, and showcase its application to parabolic equations such as the Black-Scholes equation.

The rest of the paper is organized as follows. In Section 2, we provide a detailed presentation of the DFVM for solving PDEs. In Section 3, we present numerous numerical experiments to demonstrate the effectiveness of the DFVM. Finally, we conclude with a brief summary and discussion in Section 4.

\section{The deep finite volume method (DFVM)}\label{DFVM}\

In this section, we present the deep finite volume method (DFVM) for solving the following partial diﬀerential equation
\begin{align}\label{proba}
\mathcal{L} u & = f\quad \text{in}~\Omega\\ \label{probbd}
\mathcal{B}u & = g\quad \text{on}~\partial\Omega,
\end{align}
where $\Omega$ is a bounded domain in $\mathbb{R}^d$ with boundary $\partial\Omega$, 
${\cal L}$ and ${\cal B}$ are some given interior and boundary differential operators, respectively. 

The main purpose of this section is to train a neural network $u_\theta$ to approximate the exact solution $u$ of the problem \eqref{proba}-\eqref{probbd}. Without loss of generality, we choose $u_\theta$ as a ResNet type network 
which takes the form 
\begin{equation}\label{ResNet}
    u_\theta({\bf x})=({\bf B}_{l+1} \circ {\bf B}_{l} \circ {\bf B}_{l-2} \circ \cdots \circ {\bf B}_2 \ \circ {\bf B}_0)({\bf x}), \ {\bf x} \in \mathbb{R}^{d},
\end{equation}
where each residual block is presented as
\begin{equation*}
    {\bf h}_k = {\bf B}_k ({\bf h}_{k-2}) = \sigma({\bf W}_{k} \sigma({\bf W}_{k-1}{\bf h}_{k-2}+{\bf b}_{k-1})+{\bf b}_{k})+{\bf h}_{k-2},\ k=2, 4,\cdots, l.
\end{equation*}
and 
\begin{equation*}
    {\bf h}_{0} = {\bf B}_0({\bf x}) = \sigma({\bf W}_{0}{\bf x}+{\bf b}_{0}), \ \ \ u_\theta({\bf x}) = {\bf B}_{l+1}({\bf h}_{l}) = \sigma({\bf W}_{l+1}{\bf h}_{l} + b_{l+1}).
\end{equation*}
Here, $\{{\bf h}_k \in \mathbb{R}^{m} | k=0,2, \cdots,l\}$ are the $m-$dimensional hidden state vectors, $\sigma$ is the activation function and $\boldsymbol{\theta}=\left\{{\bf W}_{l+1}, b_{l+1},{\bf W}_k, {\bf b}_k,|k=0, 1,\cdots, l\right\}$ are parameters to be trained.

Before illustrating how to train $u_\theta$, we first introduce the numerical integration formulas used for volume and boundary integrals of control volumes in different dimensions in DFVM.

\subsection{Quadratures over control volumes and their boundaries}
For simplicity, we choose a control volume (CV) in the DFVM as a cube or ball in $\mathbb{R}^d$. Precisely, given a point  ${\bf x}_0=(x^1_0,\ldots,x^d_0)\in  \mathbb{R}^d$ and a size quantity 
$h>0$, we let $V_{{\bf x_0},h}$ be the $d$-dimensional cube 
\[
V_{{\bf x_0},h} = \prod \limits_{j=1}^d [x^j_0-h, x^j_0+h],
\]
or the $d$-dimensional  ball 
\[
V_{{\bf x_0},h} =  \{{\bf x}=(x^1,\ldots,x^d)\in \mathbb{R}^d\mid \sum_{j=1}^d(x^j-x^j_0)^2\le h^2\}.
\]
Normally,  we choose $V_{{\bf x_0},h}$ as a cube when $d\le 3$ and  a ball when $d\ge 4$. Noticing that  $V_{{\bf x_0},h}$ is not necessary in the interior of $\Omega$, the actual CV is often chose as $V=V_{{\bf x_0},h}\cap \Omega$, the part of $V_{{\bf x_0},h}$ in $\Omega$. 

Next, we illustrate how to numerically calculate an integral over $V$ or $\partial V$, the boundary of $V$.  Let $Q_r(F)=\sum_{j=1}^r w_jF(g_j), r\ge 1$ be the $r$th order Gauss quadrature to calculate the integral $\int_{-1}^1 f(x)dx$, where $-1\le g_1<\ldots<g_r\le 1$ are $r$ Gauss points and $w_j,1\le j\le r$ are corresponding weights. By an affine transformation, we can use the quadrature
\[
Q_r(F, [a,b])= \frac{(b-a)}{2}\sum_{j=1}^r w_jF(g^{[a,b]}_j), \ \ g_j^{[a,b]}=\frac{a+b}{2}+\frac{b-a}{2}g_j, 1\le j\le r
\]
to calculate the 1D integral $\int_a^b F dx$. Noticing that an interval $[a,b]$ is uniquely determined by its center $c=\frac{a+b}{2}$ and half length $\hat{h}=\frac{b-a}{2}$, so sometimes, we also denote $Q_r(F, c, \hat{h})=Q_r(F, [a,b])$ and 
$g_j^{c,\hat{h}}=g_j^{[a,b]}, 1\le j\le r$.
With this notation, a Gauss quadrature on a $d$-dimensional control volume $V=V_{{\bf x_0},h}$ can be presented as
\[
Q_r(F, V_{{\bf x_0},h})=h^d\sum_{j_1,\ldots,j_d=1}^r  w_{j_1}\cdots w_{j_d} F(g_{j_1}^{x_0^1,h},\ldots, g^{x^j_0,h}_{j_i}, \ldots, g_{j_d}^{x^d_0,h}),
\]
which can be used to calculate the integral $\int_{V_{{\bf x_0},h}} F d{\bf x}$.
Next, we present quadratures for the integral on the boundary $\partial V$.
In the case $d=2$,  $\partial V$ is the union of 4 segments and thus
\begin{align}\nonumber
\int_{\partial V} F ds \approx Q_r(F,\partial V)&=Q_r(F(x^1_0+h,\cdot), x^2_0, h)+Q_r(F(x^1_0-h,\cdot),x^2_0,h))\\ \label{q2}
&+Q_r(F(\cdot, x^2_0+h), x^1_0,h)+Q_r(F(x^2_0-h,\cdot), x^1_0,h). 
\end{align}
 In the case $d=3$,  $\partial V$ is the union of 6 squares, therefore
\begin{align}\nonumber
\int_{\partial V} F ds \approx Q_r(F,\partial V) &=  Q_r(F(x^1_0+h,\cdot,\cdot), S_1)+Q_r(F(x^1_0-h,\cdot,\cdot), S_1)+Q_r(F(\cdot, x^2_0+h,\cdot), S_2) \\
&+Q_r(F(\cdot, x^2_0-h,\cdot), S_2)+Q_r(F(\cdot, \cdot, x^3_0+h), S_3)+Q_r(F(\cdot, \cdot, x^3_0-h,\cdot), S_3),
\end{align}
where $S_1=[x^2_0-h,x^2_0+h]\times [x^3_0-h,x^3_0+h]$, $S_2=[x^1_0-h,x^1_0+h]\times [x^3_0-h,x^3_0+h]$, and $S_3=[x^1_0-h,x^1_0+h]\times [x^2_0-h,x^2_0+h]$.
 
When $d>3$,  we prefer to use the (quasi) Monte-Carlo method to calculate the integrals $\int_V F d{\bf x}$ and $\int_{\partial V} F ds$. Fixing two positive integers $J_V$ and $J_{\partial V}$ which are independent of the dimension $d$, we randomly sample $J_V$ points ${\bf x}_j,j=1,\ldots, J_V$ from the interior of $V$ and $J_{\partial V}$ points ${\bf y}_j,j=1,\ldots, J_{\partial V}$ from the boundary $\partial V$. We denote $S_V:=\{{\bf x}_j\in V|j=1,\ldots, J_V\}$ and  $S_{\partial V}:=\{{\bf y}_j\in \partial V|j=1,\ldots, J_{\partial V}\}$ as  the set of training points in $V$ and $\partial V$ respectively. Then we use the quadrature in \cite{caf1998, chen2021}
\begin{eqnarray}\label{mcq1}
\int_{V} F d{\bf x} &\approx& Q_{MC}(F, V)=\frac{|V|}{J_V}\sum_{{\bf x}\in S_V} F({\bf x}),
\end{eqnarray}
and 
\begin{eqnarray}\label{mcq2}
\int_{\partial V} F ds &\approx& Q_{MC}(F, \partial V) =\frac{|\partial V|}{J_{\partial V}}\sum_{{\bf y}\in S_{\partial V}} F({\bf y}).
\end{eqnarray}

To specify \( J_V \) and \( J_{\partial V} \), we adopt the following considerations. First, given that the size of \( V \) is typically quite small (e.g., \( h \sim 10^{-5} \)), we can select \( J_V \) as a fixed constant independent of the dimension \( d \). In practice, it is common to set \( J_V = 1 \). Second, since \( |\partial V| \sim \frac{d}{h}|V| \), a suitable choice for \( J_{\partial V} \) would be around \( d \) to balance computational efficiency and accuracy. For example, when \( V \) is a cube, setting \( J_{\partial V} = 2d \) ensures there is at least one integration point per face of its boundary. Our numerical experiments consistently show that setting \( J_{\partial V} = 2d \) typically achieves sufficient accuracy.

In summary, we let the quadrature 
\begin{equation*}
Q(F,V)=\left\{\begin{array}{ll} Q_r(F,V), &d\le 3\\
Q_{MC}(F,V), &d\ge 4  
\end{array}
\right.  \ {\rm and }\ \ 
Q(F,\partial V)=\left\{\begin{array}{ll} Q_r(F,\partial V), &d\le 3\\
Q_{MC}(F,\partial V), & d\ge 4.  
\end{array}
\right. 
\end{equation*}

\subsection{Loss}
In this subsection, we illustrate how to construct the loss of the DFVM for solving PDEs.

Similar to other deep PDE solvers like PINN, the DFVM formulates its loss function in the least squares sense. Initially, we sample \( {\cal S}_{int} \), which consists of training points within the interior of \( \Omega \), and \( {\cal S}_{bdy} \), which comprises training points located on the boundary \( \partial \Omega \), according to a specified distribution (e.g., Uniform or Gaussian).

With a fixed size \( h > 0 \), a \emph{control volume} \( V_{\bf x} = V_{{\bf x},h} \cap \Omega \) is constructed around each point \( {\bf x} \in {\cal S}_{int} \). Given the random nature of \( {\cal S}_{int} \), the resulting collection of control volumes, \( \{ V_{\bf x} \ | \ {\bf x} \in \mathcal{S}_{int} \} \), does not partition the spatial domain and may include overlapping regions between adjacent volumes, depending on the size parameter \( h \). The fact that \( \{ V_{\bf x} \ | \ {\bf x} \in \mathcal{S}_{int} \} \) does not necessarily constitute a partition of the domain \( \Omega \) distinguishes DFVM from traditional finite volume methods \cite{XZ2009} 
which suffer from the {\it curse of dimensionality.}

\subsubsection{Divergence-form second-order PDEs}
We begin with  the case that $\mathcal{L}$ in \eqref{proba} is a divergence-form second-order operator 
given by
 \begin{equation}\label{div2PDE}
 \mathcal{L} u =-\nabla \cdot \left( \mathbf{A} \nabla u \right)+ {\bf b}\cdot\nabla u + cu, 
 \end{equation}
where both the matrix ${\bf A}$ and vector ${\bf b}$ are known variable-coefficients. The boundary condition(s) in \eqref{probbd} can be Dirichlet, Neumann, and Robin types.
We suppose that the coefficient matrix ${\bf A}=(a_{ij}(x):1\le i,j\le d )$ is symmetric, uniformly bounded, and positive definite in the sense that there exist positive constants $\alpha,\beta$ such that
\begin{equation*}
    \alpha \xi^t\xi \le \xi^t A(x) \xi \le \xi^t\xi, \forall x\in \Omega, \xi\in \mathbb{R}^d. 
\end{equation*}
Note that under the above properties on ${\bf A}$ and some appropriate properties on ${\bf b}$ and $c$, the corresponding PDE \eqref{proba} and \eqref{probbd} has a unique solution.

We have
\begin{align}\nonumber
    \int_{V_{\bf x}} \big( \mathcal{L} u_\theta - f {\rm \ } \big) d {\bf x}
    & =\int_{V_{\bf x}} \big( -\nabla \cdot \left( \mathbf{A} \nabla u_\theta \right) + {\bf b}\nabla u_\theta + cu_\theta - f \big) d {\bf x} \\ \label{DFVMlossA}
    & = - \int_{\partial V_{\bf x}} \left( \mathbf{A} \nabla u_\theta \right)\cdot {\bf n} ds +\int_{V_{\bf x}} \left(  {\bf b}\nabla u_\theta + cu_\theta - f \right) d {\bf x},
\end{align}
where ${\bf n}$ is the unit normal outward $V_{\bf x}$, and in the second equality, we have used the divergence theorem \cite{divergence} to transform an integral in a volume to an integral on its boundary surface.

The interior loss is then defined by 
\begin{equation}\label{DFVMintloss}
\mathcal{J}^{int}(u_\theta) = \frac{1}{\# {\cal S}_{int}} \sum_{{\bf x} \in {\cal S}_{int}}  \frac{1}{|V_{\bf x}|^2}\big| Q(\left( -\mathbf{A} \nabla u_\theta \right)\cdot {\bf n},\partial V_{\bf x})+Q({\bf b}\cdot \nabla u_\theta + cu_\theta - f , V_{\bf x})\big|^2. 
\end{equation}
We emphasize that in the loss \eqref{DFVMintloss}, no second-order derivative term is involved.

Next we define the boundary loss as 
\begin{equation}\label{DFVMbdyloss}
\mathcal{J}^{bdy}(u_\theta) = \frac{1}{\# {\cal S}_{bdy}} \sum_{{\bf x} \in {\cal S}_{bdy}}  \left|\mathcal{B} u_\theta({\bf x}) - g({\bf x})\right|^2. 
\end{equation}
The total loss function is then defined by 
\begin{equation}\label{DFVMloss}
\mathcal{J}(u_\theta) = \mathcal{J}^{int}(u_\theta)+\lambda \mathcal{J}^{bdy}(u_\theta),  
\end{equation}
where $\lambda$ represents the weight of the boundary loss term that needs to be determined.

\subsubsection{General second-order PDEs}
In this subsection, we discuss how to design the loss for  the case that $\mathcal{L}$ in \eqref{proba} is a general second-order operator 
given by
\begin{equation}\label{ndiv2PDE}
 \mathcal{L} u =-\mathbf{A} : D^2u + {\bf b}\cdot\nabla u + cu, 
 \end{equation}
 where the tensor product
\[ 
{\bf A}:D^2 v=\sum_{1\le i,j\le d} a_{ij}\partial^2_{x_ix_j} v, \forall v \in H^2(\Omega).
\]
In addition, we often assume that the coefficient tensor ${\bf A}$ satisfies the Cordes condition; that is, there exists an $\epsilon\in [0,1]$ such that
\[
\frac{|{\bf A}|^2}{(tr{\bf A})^2} \le \frac{1}{(d-1+\epsilon)},
\]
where $|{\bf A}|^2=\sum_{i,j=1}^d a^2_{ij}.$ 
Note that the Cordes condition is often necessary to ensure that the original PDE has a unique solution \cite{cordes}. 

If the coefficient matrix ${\bf A}\in [C^1(\Omega)]^{d\times d}$, then the operator ${\cal L}$ can be rewritten as the divergence form 
\begin{equation}\label{div2PDE2}
 \mathcal{L} u =-\nabla \cdot \left( \mathbf{A} \nabla u \right)+ (\nabla\cdot {\bf A}+ {\bf b})\cdot\nabla u + cu,
\end{equation}
so that we can design the loss according to the method presented in Section 2.2.1. In the case ${\bf A}\notin [C^1(\Omega)]^{d\times d}$, we do not have the above global transformation to allow us to transform the integral of all second order derivative terms on a volume to the integral of first order derivative terms on the boundary of the volume only by once. Fortunately, we can use the fact that
\[
\partial^2_{x_i,x_j}u={\rm div} {\bf v}
\]
where ${\bf v}=(v_1,\ldots, v_d)$ is a vector given by 
\[
v_k=\left\{
\begin{array}{ll}
0,   &k\not= i,\\
\frac{\partial u}{\partial x_j},& k=i,
\end{array}
\right.
\]
to transform the integral 
\[\int_{V_{{\bf x}}} \alpha_{ij} \partial^2_{x_i,x_j}u_\theta d {\bf x} \approx \alpha_{ij}({\bf x})\int_{\partial V_{{\bf x}}}  {\bf v}\cdot \vec{n} ds=\alpha_{ij}({\bf x})\int_{\partial V_{{\bf x}}}  \frac{\partial u_\theta}{\partial x_j} n_i ds,
\]
where $n_i, 1\le i\le d$ is the $i$th component of the normal vector ${\vec{n}}$.
Then the interior loss is  defined by 
\begin{equation}\label{DFVMintloss2}
\mathcal{J}^{int}(u_\theta) = \frac{1}{\# {\cal S}_{int}} \sum_{{\bf x} \in {\cal S}_{int}}  \frac{1}{|V_{{\bf x}}|^2}\big| \sum_{i,j=1}^d \alpha_{ij}({\bf x}) Q(\frac{\partial u_\theta}{\partial x_j} n_i,\partial V_{{\bf x}})+Q({\bf b}\cdot \nabla u_\theta + cu_\theta - f , V_{{\bf x}})\big|^2. 
\end{equation}

\subsubsection{High order PDEs}
If ${\cal L}$ is a differential operator of an order higher than 2, we may use some so-called mid-variables to transform \eqref{proba} to a system of second-order equations.
For instances, when ${\cal L} u=\Delta^2 u$, we introduce the mid-variable $v=\Delta u$ to convert the biharmonic equation
\begin{equation}\label{eq:bih}
    \Delta^2 u= f   
\end{equation}
into a system of two second-order equations as below
\begin{eqnarray}\label{eq:bih-dec}
    \Delta u= v   \ \ {\rm in}\ \ \Omega,  \\
    \Delta v= f   \ \ {\rm in}\ \ \Omega.
\end{eqnarray}

When ${\cal L} u =\frac{\partial u}{\partial t}+\varepsilon^{2} \Delta^{2} u-\Delta\left(u^{3}-u\right)$ is the fourth-order Cahn-Hilliard type operator, we may use the mid-variable $v=-\varepsilon^2 \Delta u+u^3-u$ to convert the Cahn-Hilliard equation
\begin{equation}\label{eq:CH}
\frac{\partial u}{\partial t}=-\varepsilon^{2} \Delta^{2} u+\Delta\left(u^{3}-u\right)+g, 
\end{equation}
into the system of second-order equations as below
\begin{align}\label{eq:CH-dec}
\varepsilon^{2} \Delta u &= -v + u^{3}-u, \\
\Delta v &= \frac{\partial u}{\partial t} - g. 
\end{align}

When ${\cal L} u =\frac{\partial u}{\partial t}-\Delta\left[u^{2}+u^{3}+\left(\left(q_{0}+\Delta\right)^{2}-\varepsilon\right) u\right] $ is the sixth-order Phase-Field Crystal operator, we may use the mid-variables $v=q_0 u+\Delta u, w=u^2+u^3-\varepsilon u+(q_0+\Delta)v $ to convert the Phase-Field Crystal equation
\begin{equation}\label{eq:PFCE}
    \frac{\partial u}{\partial t}=\Delta\left[u^{2}+u^{3}+\left(\left(q_{0}+\Delta\right)^{2}-\varepsilon\right) u\right], 
\end{equation}
where $q_{0}$ and $\varepsilon$ are constants, into the system of three second-order equations as below
\begin{align}\label{eq:PFCE-dec}
    (q_0+\Delta) u &= v , \\
    (q_0+\Delta) v &= w - u^{2} - u^{3} +\varepsilon u, \\
    \Delta w &= \frac{\partial u}{\partial t}.
\end{align}

\subsection{The adaptive trajectories sampling DFVM (ATS-DFVM)}\ 
In the previous section, we train the parameters of a neural network solution on fixed training points. In this section, we explain how to update the set of training points adaptively according to the computed approximate solution to improve the performance of a deep learning solver for PDEs.

We recall that the adaptive selection of training points is an important tool to improve the performance of a deep solver for PDEs. Along this direction, a lot of effort has been put into, see e.g. \cite{LMM2021,NGM2021,FI-PINN,ADN,DAS-PINN} for an uncompleted list of publications.
In this paper, we demonstrate how to enhance the performance of DFVM by incorporating a novel adaptive sampling technique known as ATS, which is recently developed in \cite{ATS}.
 Without loss of generality, in the following we only explain how to generate ${\cal S}^{new}_{int}$, the set of interior training points for the next training stage, adaptively  according to ${\cal S}_{int}$, the set of interior training points for the current training stage, and $u_\theta$, the neural network solution which has been trained using the points in ${\cal S}_{int}$.
To this end, we first generate a set of candidate training points.
For each ${\bf x}_i\in {\cal S}_{int}, 1\le i\le I$, where $I=\#{\cal S}_{int}$ is the cardinality of ${\cal S}_{int}$, we use the Gaussian stochastic process to generate $J$ novel points. Precisely, we let 
\begin{equation}\label{ATS-spread}
\mathbf{x}_{i, j}=\mathbf{x}_{i}+\sqrt{\Delta t} \mathcal{N}\left(0, \mathcal{I}_{d}\right), j=1, \ldots, J, 
\end{equation}
where $\Delta t>0$ is a small prescribed radius and $\mathcal{N}$ is the normalized Gauss process.
We define the next step's set of candidate training points as
\[
{\mathcal S}'_{int}=\{{\bf x}_{i,j}| 1 \leq i \leq I, 1 \leq j \leq J\}\cup {\mathcal S}_{int}.
\]
Secondly, we define a DFVM type error indicator. For a point ${\bf x}\in {\mathcal S}'_{int}$, let $V=V_{{\bf x},h}$ be a control volume defined in Section 2.1 and we define
\begin{equation}\label{ind-DFVM}
     {\rm Ind}_{V} ({\bf x}) = \big| Q(\left( -\mathbf{A} \nabla u_\theta \right)\cdot {\bf n},\partial V_{{\bf x}, h})+Q({\bf b}\cdot \nabla u_\theta + cu_\theta - f , V_{{\bf x}, h})\big|.
\end{equation}
Note that the DFVM-type error indicator described above is defined for divergence form second-order PDEs. However, it can be extended to general second-order PDEs and higher-order PDEs, similar to the approach we used for defining the loss in Section 2.2.
Finally, we construct ${\mathcal S}^{new}_{int}$ by selecting  $I$ points from  ${\mathcal S}'_{int}$, where the indicator value is largest,  to form ${\cal S}^{new}_{int}$, the set of training points in the next training stage.
That is, the next stage's  training points set ${\mathcal S}^{new}_{int}$ satisfies the following two properties : 1) $\#{\mathcal S}^{new}_{int}=\#{\mathcal S}_{int}$, 2) ${\rm Ind}_{V} ({\bf x})\ge {\rm Ind}_{V} ({\bf y})$ for all ${\bf x}\in {\mathcal S}^{new}_{int}, {\bf y}\in {\mathcal S}'_{int}\setminus {\mathcal S}^{new}_{int}$.
Note that the basic idea behind this strategy is to focus training more intensively on locations where the error indicator is relatively large. For further details on the Adaptive Training Strategy (ATS), please refer to the manuscript available at https://arxiv.org/abs/2303.15704.

\subsection{Complexity comparison: computing $\Delta u_\theta$ with DFVM vs AD}

In this subsection, we evaluate the computational complexity involved in calculating \(\Delta u_\theta\) using the DFVM method (as per scheme \eqref{disLaplace1}) and directly applying AD. We begin by analyzing the computational cost of evaluating the function value, first-order derivatives, and second-order derivatives. Following this, we provide a comparative analysis of the computational costs associated with using DFVM and AD to compute \(\Delta u_\theta\).


For simplicity, let's consider a fully connected network 
$u_{\theta}({\bf x})={\bf F}_{L} \circ {\bf F}_{L-1} \circ \cdots \circ {\bf F}_{0} ({\bf x})$, where
\begin{align*}
    {\bf h}^{0} = {\bf x}, \ \  {\bf h}^{l+1} ={\bf F}_l({\bf h}^l) = \sigma(W^l {\bf h}^l + {\bf b}^l), l=0,1,\ldots, L-1,
\end{align*}
where ${\bf x} \in \mathbb{R}^d$ and $y = u_{\theta}({\bf x})= {\bf F}_{L}({\bf h}^L)=W^L {\bf h}^L + {\bf b}^L$.  Here we let the dimension of ${\bf h}^l$ be $d_l = m$ for all $l=1,2,\ldots, L$. 

In {\it PyTorch}, the function value $y=u_\theta({\bf x)}$ is computed in the following forward propagation way: 
\[
{\bf x} \to \cdots \to {\bf h}^{l} \to {\bf h}^{l+1} \to \cdots \to y.
\]
Namely, $y$ is computed by the following iterative algorithm:
\[
{\bf h}^{l+1}={\bf F}_l({\bf h}^l), \ l=0,\ldots, L.
\]
Therefore, the evaluation of $u_\theta({\bf x})$ takes a total of $\mathcal{O}(md+Lm^2)$ multiplication operations, $\mathcal{O}(md+Lm^2)$ addition operations and $Lm$ calls of $\sigma$.
The gradient of $u_\theta$ is computed in the following backward propagation way
\begin{align*}
&\nabla_{{\bf h}^{L}}y \to \cdots \to \nabla_{{\bf h}^{l}}y \to \cdots \to \nabla_{{\bf x}}y.
\end{align*} 
That is, it is computed with the iterative scheme 
\[
 \nabla_{{\bf h}^l}y = \nabla_{{\bf h}^{l}}{\bf F}_{l} \nabla_{{\bf h}^{l+1}} y, \ l=0,\cdots,L.
\]
Therefore, the computation of the gradient of $u_\theta$ takes a total of $\mathcal{O}(md+m^2)$ multiplications , $\mathcal{O}(md+m^2)$ addition operations, and $Lm$ calls of the function $\sigma'$. The Laplacian of $y=u_\theta({\bf x})$ is practically computed by taking the trace on the fully Hessian matrix. The calculation of the Hessian matrix is obtained through the following two propagation ways
\begin{align*}
    &\nabla_{\bf x} {\bf h}^0 \to \cdots \to \nabla_{\bf x} {\bf h}^l \to \cdots \to \nabla_{\bf x} {\bf h}^L, \\
    &\nabla\nabla_{{\bf h}^L}y \to \cdots \to \nabla\nabla_{{\bf h}^{1}}y \to \cdots \to \nabla_{\bf x}\nabla y = \nabla\nabla y,
\end{align*}
The corresponding iterative schemes are 
\begin{align*}
    &\nabla {\bf h}^{l+1} = \nabla {\bf h}^{l} \nabla_{{\bf h}^{l}} {\bf F}_{l}, \\
    &\nabla \nabla_{{\bf h}^{l}} y = \nabla \nabla_{{\bf h}^{l+1}} y \nabla_{{\bf h}^{l}}{\bf F}_{l}^{T} + [(\nabla {\bf h}^{l})\nabla^2_{{\bf h}^{l}}{\bf F}_{l}] \nabla_{{\bf h}^{l+1}} y.
\end{align*}
Therefore, the computation of the Laplacian of $u_\theta$ takes a total of 
$\mathcal{O}(md^3+md^2+Lm^3d)$ multiplications, 
$\mathcal{O}(md^3+md^2+Lm^3d)$ addition operations, 
and $Lm$ calls of $\sigma''$.

The DFVM transforms the volume integral of the second-order Laplacian operator into a boundary integral involving only first-order derivatives. In actual calculations, numerical integration is used instead of the boundary integral. To achieve higher accuracy, the number of integration interpolation points used in the article is typically $2d$, corresponding to the number of sub-boundaries of a cube. Therefore, the total computational effort is $2d$ times the computation for the first-order derivatives, which is $\mathcal{O}(2md^2+2m^2d)$. Compared to directly using automatic differentiation for calculating second-order differential operators, the computational cost of DFVM is theoretically lower.

\section{Numerical experiments}\label{numerical-experiments}\

In this section, we test the performance of the DFVM. 
First, we compare the performance of the scheme \eqref{disLaplace1} and the AD by applying them to calculate $\Delta u_\theta$, the Laplacian of a neural network function $u_\theta$. With an appropriately chosen volume size ($h=10^{-5})$,  the scheme \eqref{disLaplace1} computes a very accurate approximation of $\Delta u_\theta$ by consuming far less time than that of the AD. 
Secondly, we apply the DFVM to solve variants of PDEs including the Poisson equation, the biharmonic equation, the Cahn-Hilliard equation, and the Black-Scholes equation. The numerical experiments are implemented using Python with the PyTorch library on a machine equipped with NVIDIA TITAN RTX GPUs, except for Case 3 of the Poisson equation, which is implemented on a machine equipped with an NVIDIA V100 GPU for a fair comparison.

In all our numerical experiments, unless otherwise stated, all methods of comparison share the same neural network and training points.
For all methods, we use the Adam optimizer \cite{Adam} to train the network parameters. Moreover, the activation function will be chosen as the \emph{tanh} function.
The accuracy of $u_\theta$  is indicated by the {\it $L^2$ relative error} defined by 
   \( {\rm RE}=\|u_\theta-u\|_{L^2}/\|u\|_{L^2}.\)
The random seed is set to 0.
Our code is available at \href{https://github.com/Sysuzqs/DFVM}{https://github.com/Sysuzqs/DFVM}.

\subsection{Computing the Laplacian operator based on DFVM}\label{sec:laplacian}

In this subsection, we justify our complexity analysis in Section 2.4 by numerical experiments.
To this end, we refine the scheme \eqref{disLaplace1} as follows: 
\begin{equation}\label{DFVMLaplace}
    \Delta u_{\theta}^{\text{DFVM}} ({\bf x_0})= \frac{1}{2h} \sum_{i=1}^d \left(\frac{\partial u_{\theta}}{\partial \bf{n}} ({\bf x}_0+h{\bf e}_i)-\frac{\partial u_{\theta}}{\partial \bf{n}} ({\bf x}_0-h{\bf e}_i)\right),
\end{equation}
where ${\bf e}_i$ is the unit vector whose $i$th component is $1$ and other components are zero. Note this scheme is derived from \eqref{disLaplace1} by choosing $V_{\bf x_0}$ as a cube centered at ${\bf x}_0$ and parameterized with a size quantity $h$. Additionally, a numerical quadrature point is chosen on each face of the cube. Practically, each  $\frac{\partial u_{\theta}}{\partial \bf{n}} ({\bf x}_0+h{\bf e}_i)$ will be computed by AD directly. 

We use the above scheme  to calculate  $\Delta u_\theta$ at 100 randomly chosen points in $R^d$ for variants of dimension $d$. Our to-be tested neural network function $u_\theta$ contains three residual blocks, with 128 neurons in each layer, and its activation function is chosen to be the \emph{tanh} function. We initialize the network parameters $\theta$ using a Gaussian distribution with a mean of 0 and a variance of 0.1. We use the data type {\it double} in the calculation process in our Torch implementation, which is a 64-bit precision floating point number.

We first test the accuracy of  $\Delta u_\theta^{\text{DFVM}}$ in variants of cases. Precisely, we will test the cases that the dimension $d$ varies from $2$ to $100$ and the volume size $h$ varies from $10^{-1}$ to $10^{-10}$. 
Listed in Table \ref{table:radius} are the mean absolute errors (MAEs,\cite{MAE2013}) between the values of $\Delta u_\theta$ computed by the AD directly (ground truth) and that by the DFVM scheme \eqref{DFVMLaplace}, which are computed as follows
\begin{equation*}
    \text{MAEs} = \frac{1}{100} \sum_{i=1}^{100} | \Delta u_{\theta}({\bf x}_i) - \Delta u_{\theta}^{\text{DFVM}}({\bf x}_i) |.
\end{equation*}
\begin{table}[H]
    \caption{MAEs of $\Delta u_\theta$ between the AD and the DFVM}
    \label{table:radius}
    \centering
    \begin{tabular}{c|ccccccc}
    $h$   & d=2      & d=10     & d=20     & d=40     & d=60     & d=80     & d=100    \\
    \hline
    1E-01 & 1.77E-03 & 4.67E-03 & 3.76E-03 & 3.23E-03 & 2.39E-03 & 2.43E-03 & 2.26E-03 \\
    1E-02 & 1.77E-05 & 4.69E-05 & 3.77E-05 & 3.24E-05 & 2.40E-05 & 2.43E-05 & 2.26E-05 \\
    1E-03 & 1.78E-07 & 4.69E-07 & 3.77E-07 & 3.24E-07 & 2.40E-07 & 2.43E-07 & 2.26E-07 \\
    1E-04 & 1.78E-09 & 4.69E-09 & 3.77E-09 & 3.24E-09 & 2.40E-09 & 2.43E-09 & 2.26E-09 \\
    1E-05 & 2.88E-11 & 6.86E-11 & 6.68E-11 & 9.45E-11 & 8.93E-11 & 9.35E-11 & 9.75E-11 \\
    1E-06 & 1.86E-10 & 4.51E-10 & 6.47E-10 & 8.59E-10 & 8.57E-10 & 9.75E-10 & 9.58E-10 \\
    1E-07 & 1.81E-09 & 4.80E-09 & 6.08E-09 & 8.10E-09 & 8.52E-09 & 1.06E-08 & 1.06E-08 \\
    1E-08 & 2.46E-08 & 4.48E-08 & 6.22E-08 & 7.96E-08 & 8.98E-08 & 1.03E-07 & 9.76E-08 \\
    1E-09 & 2.10E-07 & 4.99E-07 & 6.68E-07 & 7.86E-07 & 8.37E-07 & 9.41E-07 & 1.08E-06 \\
    1E-10 & 1.86E-06 & 4.41E-06 & 5.74E-06 & 8.10E-06 & 8.49E-06 & 1.07E-05 & 1.03E-05 \\
    \hline
    \end{tabular}
\end{table}
From the above table, we observe that for all dimensions,  the MAE first decreases and then increases as the radius $h$ decreases. This might be because similarly to a difference method, theoretically, the smaller the $h$, the more accurate the DFVM-calculated value approximates the exact $\Delta u_\theta$, practically, along with the decrease of the size $h$,  the accumulation error from floating-point arithmetic by the computer also increases. Therefore, the best approximation is often achieved when $h$ is neither too big nor too small. From Table \ref{table:radius}, we observe that for all dimensions, the minimum MAE is achieved when the radius  $h=10^{-5}$. 
Therefore,  unless otherwise specified, we will set $h=10^{-5}$ for all subsequent numerical experiments in this section. Remark that when $h=10^{-5}$, the MAE between the approximate value by the DFVM and the exact $\Delta u_\theta$ achieves the order of $10^{-11}$, which is sufficiently accurate.
 
Next we compare the computational cost of the AD and DFVM. Recorded in Table \ref{table:FV} are computation time by using both methods to calculate $\Delta u_\theta $ 10,000 times over 100 randomly chosen testing points. 
\begin{table}[H]
    \caption{Computation time by the DFVM and the AD}
    \label{table:FV}
    \centering
    \begin{tabular}{lccc}
    d   & MAE      & AD time (s) & DFVM time (s) \\
    \hline
    2   & 2.88E-11 & 36   & 12   \\
    4   & 3.13E-11 & 60   & 18   \\
    8   & 4.23E-11 & 105  & 24   \\
    10  & 6.86E-11 & 127  & 25   \\
    20  & 6.68E-11 & 265  & 47   \\
    40  & 9.45E-11 & 590  & 109   \\
    50  & 9.83E-11 & 614  & 118   \\
    60  & 8.93E-11 & 834  & 170   \\
    80  & 9.35E-11 & 1112 & 187   \\
    100 & 9.75E-11 & 1344 & 230  \\
    \hline
    \end{tabular}
\end{table}
From this table, we find that for all dimensional cases, the consumed computing time of the DFVM is far less than that of the AD. In particular, for the cases $d=80,100$, the computation time by the DFVM is almost only $1/6$ of that by the AD, while the MAE achieves 
$10^{-11}$ which means that the approximate value of $\Delta u_\theta^{\text{DFVM}}$ calculated by the DFVM is very accurate.

\subsection{Solving PDEs with the DFVM}
\subsubsection{The Poisson equation} 
We consider the Poisson equation with Dirichlet boundary condition 
which has the following form
\begin{equation}\label{eq:Poisson}
        -\Delta u = f \ {\rm in \ } \Omega, \ \ \ 
        u = g \ {\rm on \ } \partial \Omega,
\end{equation}
where $\Omega, f, g$ will be specified in the following four cases.

{\bf Case 1}
In the first case, we let $\Omega=(0,1)^{2} \subset \mathbb{R}^{2}$, $f \equiv-2 $ in $ \Omega$ , and  $g\left(x_{1}, 0\right)=g\left(x_{1}, 1\right)=x_{1}^{2} $ for $ 0 \leq x_{1} \leq \frac{1}{2}, g\left(x_{1}, 0\right)=g\left(x_{1}, 1\right)=\left(x_{1}-1\right)^{2}$  for  $\frac{1}{2} \leq   x_{1} \leq 1 $, and  $g\left(0, x_{2}\right)=g\left(1, x_{2}\right)=0$  for  $0 \leq x_{2} \leq 1 $ on $ \partial \Omega$. In this case, the  problem \eqref{eq:Poisson} admits the solution
\begin{equation*}
    u^{*}(\mathbf{x})=u^{*}\left(x_{1}, x_{2}\right)=\left\{\begin{array}{ll}
    x_{1}^{2}, & 0 \leq x_{1} \leq \frac{1}{2}, \\
    \left(x_{1}-1\right)^{2}, & \frac{1}{2} \leq x_{1} \leq 1,
    \end{array}\right.
\end{equation*} 
which is continuous but nonsmooth.

\begin{figure}[H]
    \centering
    \subfigure[Exact solution]{
        \begin{minipage}[t]{0.22\linewidth}
            \centering
            \includegraphics[width=0.95\linewidth]{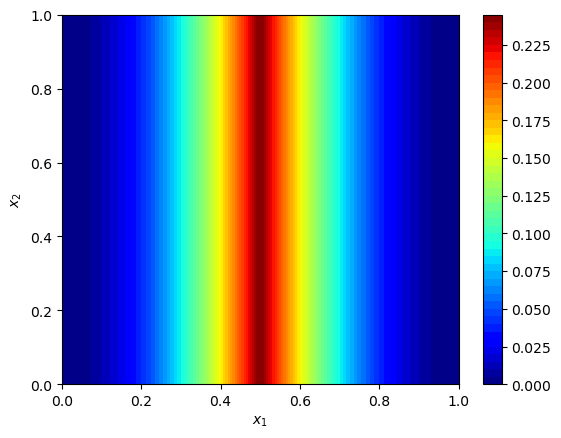}
        \end{minipage}
    }%
    \subfigure[DFVM solution]{
        \begin{minipage}[t]{0.22\linewidth}
            \centering
            \includegraphics[width=0.95\linewidth]{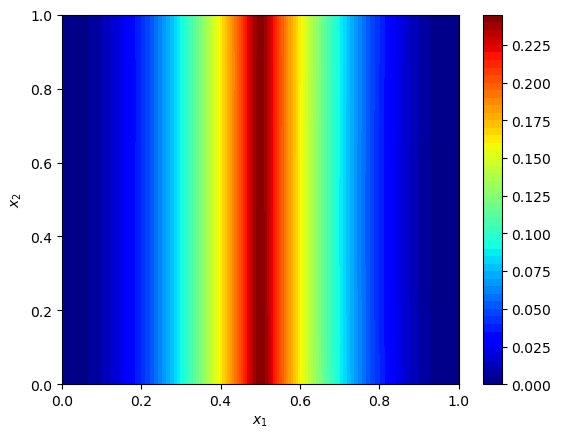}
        \end{minipage}
    }%
    \subfigure[PINN solution]{
        \begin{minipage}[t]{0.22\linewidth}
            \centering
            \includegraphics[width=0.95\linewidth]{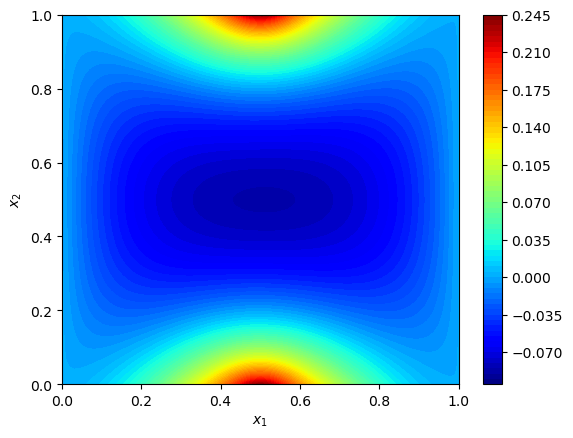}
        \end{minipage}
    }%
    \subfigure[DRM solution]{
        \begin{minipage}[t]{0.22\linewidth}
            \centering
            \includegraphics[width=0.95\linewidth]{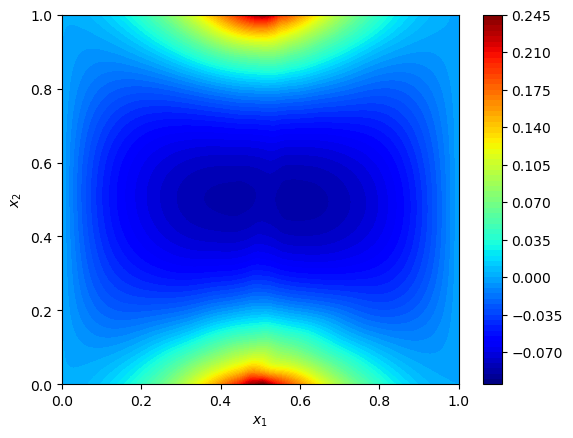}
        \end{minipage}
    }%
    \centering
    \caption{The exact solution and approximate solution by different learning methods.}
    \label{fig:singularity}
\end{figure}

\begin{figure}[H]
    \centering
    \subfigure[Loss w.r.t Iteration]{
        \begin{minipage}[t]{0.36\linewidth}
            \centering
            \includegraphics[width=0.95\linewidth]{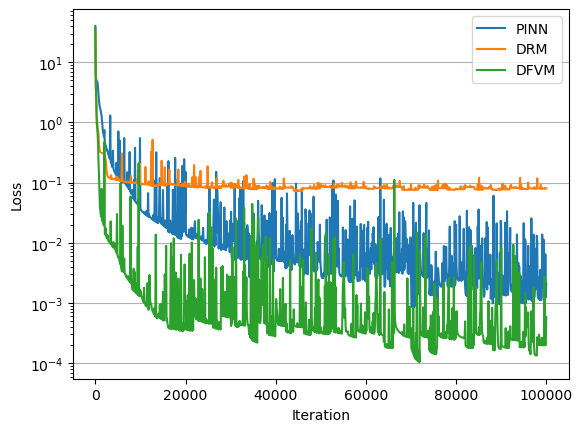}
        \end{minipage}
    }%
    \subfigure[Related error w.r.t Iteration]{
        \begin{minipage}[t]{0.36\linewidth}
            \centering
            \includegraphics[width=0.95\linewidth]{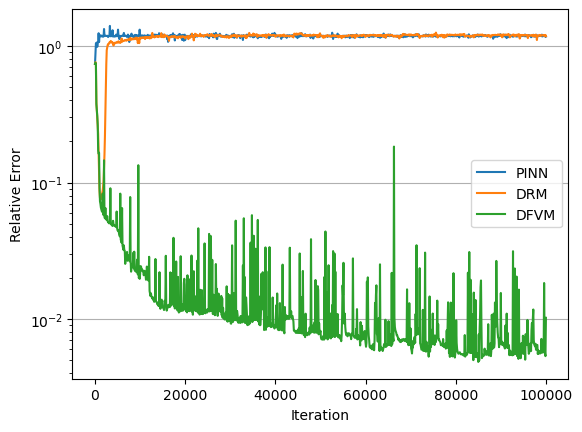}
        \end{minipage}
    }%
    \centering
    \caption{Training curves of three methods for \eqref{eq:Poisson}.}
    \label{fig:SINGcurve}
\end{figure}

We choose $u_\theta$ as a fully-connected network that has 6 hidden layers, with 40 neurons per hidden layer. To train $u_\theta$, the number of internal points is set to be 10,000, and the number of boundary points is set to be 400. The weight of the boundary loss term is set to 1000. For the DFVM, we set the control volume radius h to $10^{-3}$, $J_V=1$, $J_{\partial V}=4$. The optimizer used in each method is Adam with a learning rate of 0.001. 

Listed in Table \ref{table:sing} are the $L^2$ relative errors (REs) between the exact solution and 
the approximate solution obtained from training 100,000 iterations by the PINN, DRM, and DFVM respectively. We find that for this case in which the exact solution is nonsmooth, the DFVM 
achieves much better accuracy with much less computation time than that by the PINN. To deepen impression, we depict the exact solution and the approximate solutions by the DFVM, the PINN, and the DRM in Figure \ref{fig:singularity}. We find that the DFVM solution approximates the exact solution very closely, but the other two methods do not. 
Figure \ref{fig:SINGcurve} displays the convergence history of the loss function and relative error over 100,000 iterations. The plot illustrates that while the loss curve reaches convergence early in the training process, the relative error continues to decrease steadily.
\begin{table}[H]
    \caption{REs and computation time of case 1 after 100,000 iterations}
    \label{table:sing}
    \centering
    \begin{tabular}{ccc}
         & RE       & Time (s) \\
    \hline
    PINN & 9.32E-01 & 1033    \\
    DRM  & 8.88E-01 & 615     \\
    DFVM & 4.56E-03 & 600     \\
    \hline
    \end{tabular}
\end{table}

{\bf Case 2 }\ 
In the second case, we let  $\Omega = (0, 1)^d, f({\bf x}) = \displaystyle{\frac{1}{d} \left(\sin\left(\sum_{i=1}^{d} \frac{1}{d} x_i\right) - 2\right) }, g({\bf x}) = \displaystyle{ \left(\sum_{i=1}^{d} \frac{1}{d} x_i\right)^2 + \sin\left(\sum_{i=1}^{d} \frac{1}{d} x_i\right) }$. In this case, the solution of \eqref{eq:Poisson} is
\begin{equation*}
	 u^\ast({\bf x}) = \left(\sum_{i=1}^{d} \frac{1}{d} x_i\right)^2 + \sin\left(\sum_{i=1}^{d} \frac{1}{d} x_i\right).
\end{equation*}

In this case, the solution is sufficiently smooth. Our purpose is to test the dynamic change of the performance of the DFVM along with the dimension $d$. The network used in this case consists of 3 ResNet blocks, with each layer containing 128 neurons. Precisely, we will use the DFVM to solve \eqref{eq:Poisson} for the dimensions $d=2,4,10,20,40,60,80,100,200$. For comparison, we will compute corresponding solutions by the PINN. We set the size of the interior training points to be $2000$, and the size of the boundary training points to be $100*d$. The weight of the boundary loss term is set to 1000.

\begin{table}[H]
    \caption{REs and computation time of case 2 after 20,000 iterations}
    \label{table:PE}
    \centering
    \begin{tabular}{l|cc|cc}
        & \multicolumn{2}{c|}{DFVM} & \multicolumn{2}{c}{PINN} \\
    \hline
    d   & RE       & Time (s)       & RE       & Time (s)       \\
    \hline
    2   & 1.76E-04 & 156           & 1.73E-04 & 321           \\
    4   & 2.33E-04 & 190           & 2.44E-04 & 465           \\
    8   & 4.60E-04 & 292           & 4.71E-04 & 934           \\
    10  & 6.77E-04 & 367           & 7.23E-04 & 984           \\
    20  & 1.91E-03 & 703           & 2.22E-03 & 2173          \\
    40  & 3.40E-03 & 1424          & 4.19E-03 & 3891          \\
    50  & 5.50E-03 & 1787          & 5.00E-03 & 4450          \\
    60  & 5.52E-03 & 2163          & 5.66E-03 & 5213          \\
    80  & 6.89E-03 & 3034          & 6.05E-03 & 5856          \\
    100 & 8.02E-03 & 3899          & 7.23E-03 & 6906          \\
    200 & 8.82E-03 & 8203          & 9.11E-03 & 12552         \\
    \hline
    \end{tabular}
\end{table}
Listed in Table \eqref{table:PE} are the relative errors and computation time obtained from training $u_\theta$ 20000 iterations by the DFVM and the PINN. We observe that for all dimensional cases, the relative errors by both methods are of the same order, and the computation time of the DFVM is far less than that of the PINN.




{\bf Case 3} In this case, we let $\Omega=[-1,1]^{10}$ and the functions $f,g$ are chosen so that \eqref{eq:Poisson} admits the exact solution
\begin{equation}\label{exactsolution}
    u(\mathbf{x}) = e^{-10\|\mathbf{x}\|_2^2} .
\end{equation}
Note that for this example, the solution $u=1$ at the origin and it decays very rapidly to zero as 
the location moves from the origin to the boundary of $\Omega$.

We choose $u_\theta$ as a fully connected neural network with 7 hidden layers and  20 neurons per layer, and we choose the activation function to be the \emph{tanh} function. 
We will use three methods: the DFVM, the ATS-DFVM, and the ATS-PINN to solve the equation \eqref{eq:Poisson} in this case. 
Note that in the initial stage of both the ATS-DFVM and the ATS-PINN, we randomly select 2000 interior points and 1000 boundary points to train 3000 iterations with the PINN method. After this initial stage, we resample the training points 10 times with  ATS strategies, and after each resampling, we train the neural network 3000 iterations. The optimizer used in each method is Adam with a learning rate of 0.0001.

\begin{figure}[H]
    \centering
    \subfigure[1st resampling]{
        \begin{minipage}[t]{0.22\linewidth}
            \centering
            \includegraphics[width=0.95\linewidth]{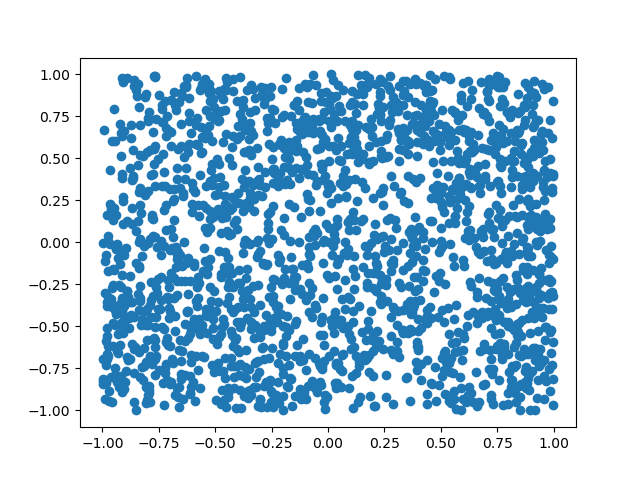}
        \end{minipage}
    }%
    \subfigure[4th resampling]{
        \begin{minipage}[t]{0.22\linewidth}
            \centering
            \includegraphics[width=0.95\linewidth]{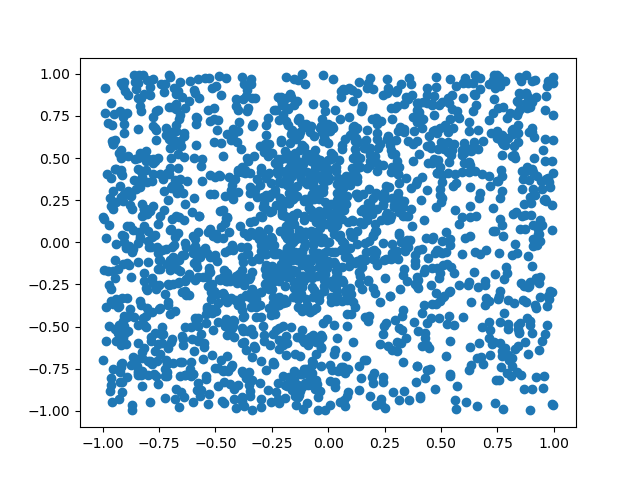}
        \end{minipage}
    }%
    \subfigure[7th resampling]{
        \begin{minipage}[t]{0.22\linewidth}
            \centering
            \includegraphics[width=0.95\linewidth]{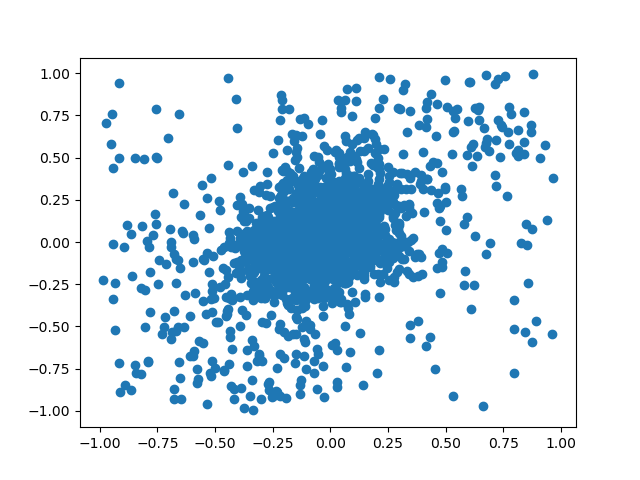}
        \end{minipage}
    }%
    \subfigure[10th resampling]{
        \begin{minipage}[t]{0.22\linewidth}
            \centering
            \includegraphics[width=0.95\linewidth]{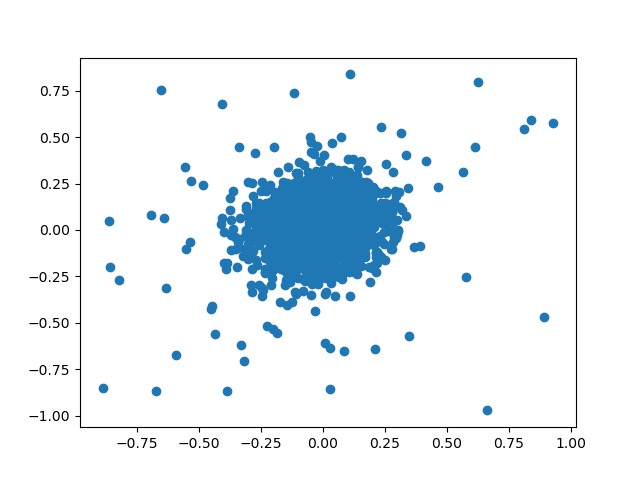}
        \end{minipage}
    }%
    \centering
    \caption{ATS-DFVM training points in the $(x_1,x_2)$ plane. For display purpose, images only show the slices of $x_3 = \cdots = x_{10} = 0$.}
    \label{fig:ATSdist}
\end{figure}

We compute the  DFVM, ATS-DFVM, and ATS-PINN solution five times (with random seeds from 0 to 4). To measure the quality of the approximation, we generate 10000 test points around the origin (in $[-0.1, 0.1]^{10}$). The average relative errors and training times of five experiments are listed in Table \ref{table:peak}. 
\begin{table}[H]
    \caption{REs and computation time of case 3 after 33,000 iterations.}
    \label{table:peak}
    \centering
    \begin{tabular}{ccc}
         & RE       & Time (s)   \\
    \hline
    DFVM      & 1.00E+00 & 465  \\
    ATS-DFVM  & 5.62E-02 & 605  \\
    ATS-PINN  & 3.82E-02 & 2323 \\
    \hline
    \end{tabular}
\end{table}
From this table, we find that the DFVM solution 
is not a good approximation of the true solution. Actually, by checking our numerical results, we find that the DFVM solution is almost zero in the whole domain. One reason for this phenomenon might be that at almost all training points, which are sampled randomly according to the uniform distribution and thus away from the origin point, the true solution is very close to zero, therefore $u_\theta$, which is trained using these points, will be also almost zero. However, by the ATS techniques, the training points will gradually close to the origin point, see Figure \ref{fig:ATSdist}. 
Consequently, we finally obtain a nice approximation of the true solution by the ATS-DFVM. From Table \ref{table:peak}, we find that the ATS techniques can also improve the approximation accuracy of the PINN. However, the computational cost of the ATS-PINN is far more than that of the ATS-DFVM. 

{\bf Case 4}
In this case, we consider a sector-shape domain $\Omega=\{(r,\theta): 0\leq r\leq 1,0\leq \theta\leq \frac{\pi}{6}\}$ and set $f=0$,
and choose $g$ in \eqref{eq:Poisson} to fit the exact solution $u(r,\theta)=r^{\frac{2}{3}}\sin(\frac{2}{3}\theta)$.
We choose $u_\theta$ as a fully connected neural network with 7 hidden layers and  20 neurons per layer, and we choose the activation function to be the \emph{tanh} function. The number of internal points is set to be 1500, and the number of boundary points is set to be 100. The weight of the boundary loss term is set to 1000. For the DFVM, we set the control volume radius $h$ to $10^{-3}$, $J_V=1$, $J_{\partial V}=4$. The optimizer used in each method is Adam with a learning rate of 0.001. 

\begin{table}[H]
    \centering
    \caption{REs and computation time of case 4 after 20,000 iterations}
    \label{table:sector}
    \centering
    \begin{tabular}{ccc}
         & RE       & Time (s) \\
    \hline
    PINN & 6.03E-03 & 243    \\
    DRM  & 7.98E-02 & 120     \\
    WAN  & 5.79E-01 & 819     \\
    DFVM & 6.02E-03 & 121     \\
    \hline
    \end{tabular}
\end{table}

Listed in Table \ref{table:sector} are the L$^2$ norm relative errors and training times for this problem. We observe that DFVM achieves the highest level of computational accuracy while maintaining high computational efficiency, demonstrating its ability to solve problems in complex domains.

\subsubsection{High order PDEs}\ 

{\bf The biharmonic equation}  We consider 
\begin{equation}\label{bih}
\Delta^2 u= f   \ \ {\rm in}\ \ \Omega    
\end{equation}
with the boundary conditions
 \(u({\bf x}) = \sum_{k=1}^{2} \text{sin}\left(\frac{\pi}{2} x_k\right),  \ \ 
\frac{\partial u}{\partial {\bf n}}({\bf x})= 0 \ {\rm on} \ \ \partial\Omega, 
\)
where $\Omega=[-1,1]^2$, and $\displaystyle{ f({\bf x})=\frac{\pi^4}{16} \sum_{k=1}^{2} \text{sin}\left(\frac{\pi}{2} x_k \right)}$. In this case, the exact solution is
\begin{equation}\label{eq:bih-solution}
    u({\bf x}) = \sum_{k=1}^{2} \text{sin}(\frac{\pi}{2} x_k).
\end{equation}

For this example, we choose $u_\theta$ as a full-connected network which has 4 hidden layers, with 40 neurons per hidden layer. To train $u_\theta$, the number of internal points is set to be 10,000, and the number of boundary points is set to be 800. The weight of the boundary loss term is set to 1000. For the DFVM, we set the control volume radius h to be 1e-3, $J_V=1$, $J_{\partial V}=4$.  The optimizer used in each method is Adam with a learning rate of 0.0001. 

\begin{table}[H]
    \caption{REs and computation time after 50,000 iterations for \eqref{bih}.}
    \label{table:bihar}
    \centering
    \begin{tabular}{ccc}
         & RE       & Time (s) \\
    \hline
    PINN & 3.27E-04 & 3647   \\
    DFVM & 1.21E-04 & 858    \\
    \hline
    \end{tabular}
\end{table}

{\bf The  Cahn-Hilliard equation}
We consider 
\begin{eqnarray}\label{CH}
\frac{\partial u}{\partial t} =-\varepsilon^{2} \Delta^{2} u+\Delta\left(u^{3}-u\right),  \ \ &{\rm in}\ \ \Omega\times [0,T],   \\
\partial_{\vec{n}} u = \partial_{\vec{n}} (-\varepsilon^2\Delta u + u^3 - u)=0,\ \ &{\rm on}\ \partial\Omega\times [0,T],\\
 u(\cdot,0) =u_0(\cdot),  \ \ \    &{\rm in}\ \Omega,
\end{eqnarray}
where $\Omega=[-1,1]^2$, $T=0.1$, $\varepsilon=0.1$, and 
\begin{equation*}
    u_0= \tanh ( \frac {1}{\sqrt {2} \varepsilon} \min \{\sqrt {(x+0.3)^ {2}+y^ {2}}-0.3,\sqrt {(x-0.3)^ {2}+y^ {2}}-0.25 \} ). 
\end{equation*}

For this example, we choose $u_\theta$ as a full-connected network that has 4 hidden layers, with 40 neurons per hidden layer. To train $u_\theta$, the number of internal points is set to 10000, the number of boundary points is set to be 800, and the number of initial points is set to be 10000. The weight of the boundary loss term is set to 1000. For the DFVM, we set the control volume radius $h$ to be $1^{-3}$, $J_V=1$, $J_{\partial V}=4$.  The optimizer used in each method is Adam with a learning rate of 0.001. 

\begin{table}[H]
    \caption{REs and computation time after 50,000 iterations for \eqref{CH}.}
    \label{table:C-H}
    \centering
    \begin{tabular}{ccc}
         & RE       & MAE \\
    \hline
    PINN & 1.09E-01 & 1.77E-01   \\
    DFVM & 7.45E-02 & 1.42E-01   \\
    \hline
    \end{tabular}
\end{table}

\begin{figure}[htbp]
    \centering
    \subfigure[DFVM solution]{
        \begin{minipage}[t]{0.7\linewidth}
            \centering
            \includegraphics[width=0.95\linewidth]{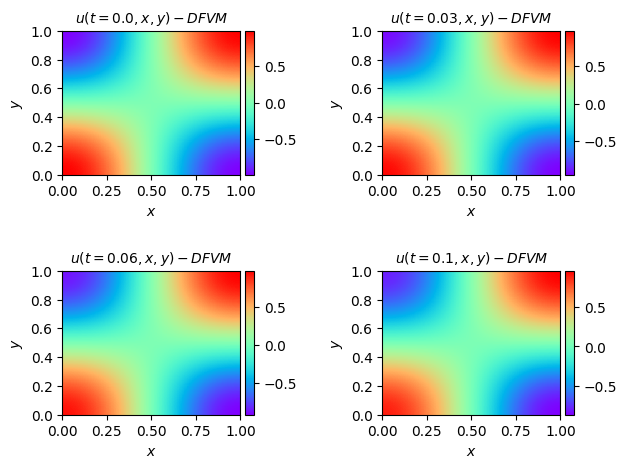}
        \end{minipage}
    }%
    \centering
    \caption{Results of C-H equation.}
    \label{fig:C-H}
\end{figure}

In this example, both DFVM and PINN demonstrate suboptimal accuracy. Recent advancements in operator-based learning methods, such as Neural Networks for learning Mean Curvature Flow (NNMCF, \cite{nnmcf}), have shown promising results for phase-field problems. Additionally, strategies like adaptive PINN \cite{APINNs} and sequential methods \cite{seqPINNs} have been successfully applied to enhance PINN performance by modifying training strategies. Importantly, these strategies can also be adapted for use with DFVM. Therefore, in this context, we specifically compare the performance of DFVM and PINN in solving the Cahn-Hilliard equation over a short time interval ($T=0.1$).




\subsubsection{Black-Scholes Equation}
In this example, We consider the well-known Black-Scholes equation below
\begin{equation}\label{eq:BSeq}
\begin{cases}
\displaystyle{
    \frac{\partial u}{\partial t}(t, x)=-\frac{1}{2} \operatorname{Tr}\left[0.16 \operatorname{diag}\left(x^{2}\right) \operatorname{Hess}_{x} u(t, x)\right]+0.05(u(t, x)-(\nabla u(t, x), x)), }\\
    u(T, x)=\|x\|^{2},
\end{cases}
\end{equation}
in $\Omega \times [0,T]$, 
which admits an exact solution
\begin{equation*}
    u(x, t)=\exp \left(\left(0.05+0.4^{2}\right)(T-t)\right)\|x\|^{2}.
\end{equation*}

The equation \eqref{eq:BSeq} has been discussed in \cite{FBSNN}, but here we discuss an alternative formulation of the same equation. 
We take $\Omega \times [0,T] =[0,2]^2 \times [0, 0.01]$ and obtain
\begin{equation}\label{eq:divBSeq}
    u_{t}=-0.08 \text{div}\begin{pmatrix}
        \displaystyle{ x_{1}^{2} \frac{\partial u}{\partial x_{1}} }\\
        \displaystyle{ x_{2}^{2} \frac{\partial u}{\partial x_{2}} }
    \end{pmatrix}
    +0.05 u+(0.16-0.05)(\nabla u, \mathbf{x}).
\end{equation}

In this test, the number of internal points is set to 1000, the number of initial points is set to 1000, and the weight of the boundary loss term is set to 1. For the DFVM, we set the control volume radius h to be 1e-5, $J_V=1$, $J_{\partial V}=4$.  The optimizer used in each method is Adam. The learning rate is set to 0.01, and it decays by 0.95 every 50 steps. The neural network is chosen as the ResNet which has 3 blocks and 64 neurons per layer.

\begin{table}[h!]
    \centering
    \caption{Errors and computation time after 20,000 iterations for \eqref{eq:BSeq}.}
    \begin{tabular}{ c c c c } 
     \toprule
        Method & RE     & RE$_{0}$ & Time (s) \\ [0.5ex] 
     \hline
        DFVM   & 1.07E-03 & 1.66E-03 & 481 \\ [0.5ex]
        PINN   & 5.46E-03 & 9.40E-03 & 613 \\ [0.5ex]
     \hline
    \end{tabular}
    \label{table:BSequation}
\end{table}

\begin{figure}[htbp]
    \centering
    \subfigure[RE w.r.t. iteration steps]{
        \begin{minipage}[t]{0.6\linewidth}
            \centering
            \includegraphics[width=0.95\linewidth]{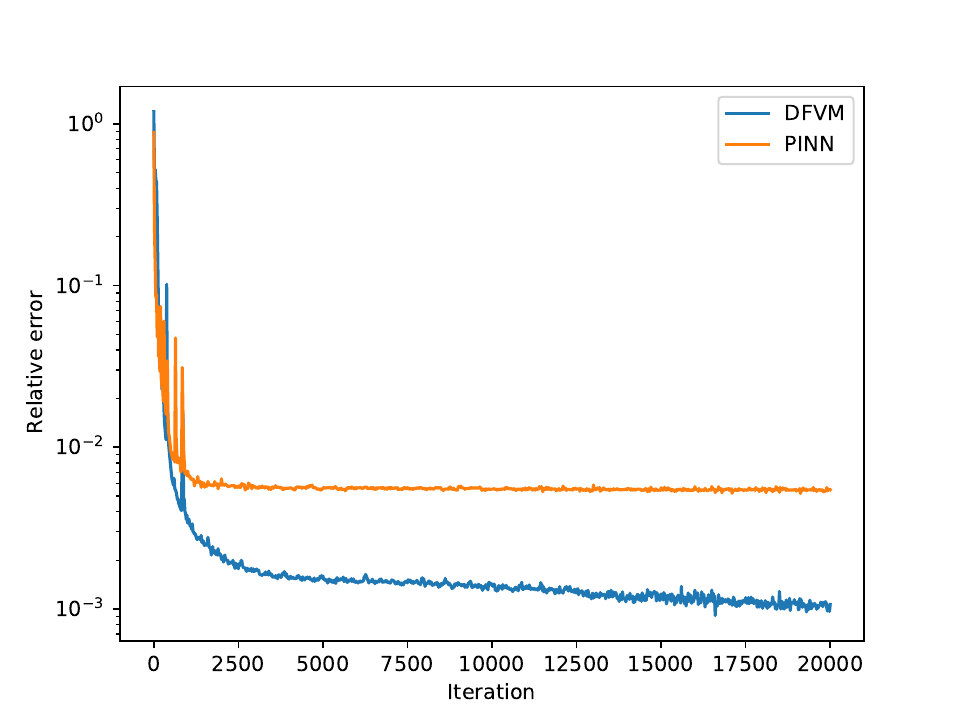}
        \end{minipage}
    }%
    \centering
    \caption{ Results of Black-Scholes Equation. }
    \label{fig:BSeq}
\end{figure}

Fig \ref{fig:BSeq} illustrates the dynamic changes of the relative errors for both DFVM and PINN methods as the iteration steps and time increase. The specific numerical values are listed in Table \ref{table:BSequation}, and we also evaluate the start point's relative error of the two methods, which is defined by $\displaystyle{ {\rm RE}_{0}=\frac{\|u_\theta(0,x)-u(0,x)\|_{L^2}}{\|u(0,x)\|_{L^2}} }$. We observe that DFVM outperforms PINN in terms of both accuracy and computational efficiency. We conclude that DFVM also performs well on parabolic PDEs.

\section{Concluding remarks} 
In the paper, we propose the DFVM, a hybrid approach merging traditional finite volume methods with deep learning techniques, which formulates the loss function based on the conservation-type weak form of PDEs rather than the traditional strong form. This approach enables DFVM to achieve higher accuracy and competitiveness in handling problems with singularities compared to deep learning methods based on the strong form. Unlike other weak forms that require surrogate test functions, DFVM directly captures local solution features, offering significant advantages in computational efficiency, accuracy, and applicability.

Moving forward, we aim to advance DFVM by integrating insights from classical numerical analysis into deeper applications. Future developments will focus on enhancing the efficiency and stability of DFVM for solving hyperbolic equations and long-time evolution problems, such as conservation laws and Navier-Stokes equations.

\section*{Acknowledgments}
The research was partially supported by the National Natural Science Foundation of China under grants 92370113 and 12071496, by the Guangdong Provincial Natural Science Foundation under the grant 2023A1515012097.

\end{document}